\documentclass[letterpaper]{amsart}

\usepackage[T1]{fontenc}
\usepackage{lmodern}
\usepackage{amssymb}
\usepackage{stmaryrd}
\usepackage{mathtools}
\usepackage{tikz-cd}
\usetikzlibrary{positioning}
\makeatletter
\tikzset{base/.code =\tikz@lib@place@handle@{#1}{base}{-1}{0}{base}{1}}
\makeatother
\usepackage[pdfusetitle]{hyperref}
\usepackage{cleveref}
\def\doi#1{\href{http://dx.doi.org/#1}{doi:#1}}

\theoremstyle{plain}
\newtheorem{theorem}{Theorem}[section]
\newtheorem{proposition}[theorem]{Proposition}
\newtheorem{corollary}[theorem]{Corollary}

\theoremstyle{remark}

\numberwithin{equation}{section}

\def\cf{\emph{cf.}}

\let\newterm\emph

\DeclareMathAlphabet\mathbfit{OML}{cmm}{b}{it}

\def\gg{\mathbfit{g}}
\def\kk{\mathbfit{k}}
\def\ll{\mathbfit{l}}
\def\bb#1{\llbracket#1\rrbracket}

\def\timesG{\mathop{\times_{G}}}
\def\timesL{\mathop{\times_{L}}}
\def\timestau{\mathop{\times_{\tau}}}
\def\timestauL{\mathop{\times_{\tau_{L}}}}
\def\timestauBL{\mathop{\times_{\tau^{BL}}}}
\DeclareMathOperator{\id}{id}

\begin{document}

\title[Fibrations of classifying spaces]{Fibrations of classifying spaces\\in the simplicial setting}
\author{Matthias Franz}
\thanks{The author was supported by an NSERC Discovery Grant.}
\address{Department of Mathematics, University of Western Ontario,
      London, Ont.\ N6A\;5B7, Canada}
\email{mfranz@uwo.ca}

% \keywords{Simplicial group, classifying space, homotopy fibration sequence,
%   twisted Cartesian product, Kan loop group, semidirect product}
% 55R10  Fiber bundles
% 55R35  Classifying spaces of groups and $H$-spaces
% 55R37  Maps between classifying spaces
% 55U10  Simplicial sets and complexes

\begin{abstract}
  In this note we show that in the simplicial setting,
  the classifying space construction converts short exact sequences of groups
  not just to homotopy fibrations, but in fact to fibre bundles.
\end{abstract}

\maketitle

\section{Introduction}

To every topological group~\(G\) one case associate its classifying space~\(BG\).
The assignment is functorial if one uses, say, Milnor's construction,
\cf~\cite[\S 14.4.3]{tomDieck:2008}.
It is well-known that if
\begin{equation}
  1 \longrightarrow K \longrightarrow G \longrightarrow L \longrightarrow 1
\end{equation}
is a short exact sequence of topological groups,
then
\begin{equation}
  BK \longrightarrow BG \longrightarrow BL
\end{equation}
is a homotopy fibration sequence.

Let us recall the short proof, \cf~\cite[pp.~348]{tomDieck:2008}:
Consider the product~\(EG\times EL\), which is another model for~\(EG\),
and its canonical projection to~\(EL\). Dividing out the actions gives a bundle
\(EG \timesG EL \longrightarrow BL\) with fibre~\(EG/K\) since
\begin{equation}
  EG \times_{G} EL = \bigl((EG \times EL)/K\bigr)/L = (EG/K)\timesL EL
\end{equation}
is the Borel construction of the \(L\)-space~\(EG/K\).
The inclusion~\(EK\to EG\) and the diagonal map~\(EG\to EG\times EL\)
give homotopy equivalences~\(BK\to EG/K\) and \(BG\to EG \timesG EL\)
which combine to the homotopy commutative diagram
\begin{equation}
  \label{eq:homfib}
  \begin{tikzcd}
    BK  \arrow{d} \arrow{r} & BG \arrow[end anchor={[xshift=-2.7pt]north}]{d} \arrow{r} & BL \arrow{d} \\
    EG/K \arrow{r} & | [base=-3.7pt] | EG \timesG EL \arrow{r} & BL \mathrlap{.}
  \end{tikzcd}
\end{equation}

While the right square commutes, the left one does so
only up to homotopy because
the image of the composition~\(EK \to EG\to EL\) is \(E1\),
which (using Milnor's construction) is an infinite-dimensional simplex
instead of a point.
If one works simplicially, however, then \(E1\) is a point,
and \eqref{eq:homfib} commutes on the nose.
Going further, one might ask whether \(BG\to BL\) is itself
a fibre bundle in the simplicial setting,
that is, a twisted Cartesian product.
After all, both~\(BG\) and any simplicial fibre bundle
with fibre~\(BK\) and base~\(BL\) are isomorphic to~\(BK\times BL\)
as graded sets. This turns out to be true:

\begin{theorem}
  \label{thm:main-intro}
  Given a short exact sequence of simplicial groups
  \begin{equation*}
    1 \longrightarrow K \longrightarrow G \longrightarrow L \longrightarrow 1,
  \end{equation*}
  the Kan loop group~\(\Omega BL\) acts on~\(BK\).
  The simplicial set~\(BG\) is isomorphic to a twisted Cartesian product
  with fibre~\(BK\), base~\(BL\) and structure group~\(\Omega BL\).
\end{theorem}

The proof appears in \Cref{sec:general},
after a review of background material in \Cref{sec:review}.
If \(G\) is the semidirect product of~\(K\) and~\(L\),
then the structure group of the bundle reduces to~\(L\) and all constructions simplify.
This is discussed in \Cref{sec:semidirect}.

\section{Review of simplicial constructions}
\label{sec:review}

From now on, we work in the simplicial category.
So `space' will mean `simplicial set' and `group' `simplicial group',
\cf~\cite[\S\S 1,~17]{May:1968}.
For the convenience of the reader we review several constructions
from~\cite[\S\S 18,~21,~26]{May:1968} that we need later on.

Note that in our notation the face and degeneracy maps of a simplicial group
bind stronger than the group multiplication. For example,
\(\partial_{0}g h\) means \((\partial_{0}g)h\).

\subsection{Fibre bundles}
\label{sec:bundles}

Let \(B\) be a space and \(G\) a group. A \newterm{twisting function} is a
collection of maps~\(\tau\colon B_{n}\to G_{n-1}\) (\(n\ge1\))
such that
\begin{align}
  \partial_{0}\tau(b) &= \tau(\partial_{0}b)^{-1}\,\tau(\partial_{1}b), \\
  \partial_{i}\tau(b) &= \tau(\partial_{i+1}b) & (1\le i\le n-1), \\
  s_{i}\tau(b) &= \tau(s_{i+1}b) & (0\le i\le n-1), \\
  \tau(s_{0}b) &= 1_{n}.
\end{align}

For any (left) \(G\)-space~\(F\) we can form the \newterm{fibre bundle}
(or `twisted Cartesian product')~\(F\timestau B\). As a graded set,
it is the Cartesian product of~\(F\) and~\(B\), but with the zeroeth
face map twisted by~\(\tau\). More precisely,
\begin{align}
  \partial_{0}(f,b) &= \bigl(\tau(b)\,\partial_{0}f,\partial_{0}b\bigr), \\
  \partial_{i}(f,b) &= (\partial_{i}f,\partial_{i}b) & (1\le i\le n)\mathrlap{,} \\
  s_{i}(f,b) &= (s_{i}f,s_{i}b) & (0\le i\le n)  
\end{align}
for~\((f,b)\in(F\timestau B)_{n}=F_{n}\times B_{n}\).

A \newterm{pseudo-cross section} of the canonical projection~\(\pi\colon F\timestau B\to B\)
is a map \(\sigma\colon B\to F\timestau B\) such that
\begin{align}
  \pi\sigma &= \id_{B}, \\
  \partial_{i}\sigma &= \sigma\partial_{i} & (1\le i\le n), \\
  s_{i}\sigma &= \sigma s_{i} & (0\le i\le n).
\end{align}
If additionally \(\partial_{0}\sigma=\sigma\partial_{0}\), then \(\sigma\)
is a simplicial map, thus an honest cross section.

\subsection{Classifying spaces}

Let \(G\) be a group. The \newterm{classifying space}~\(BG\) of~\(G\) is defined
by~\((BG)_{0}=\{[\,]\}\), \((BG)_{n}=G_{n-1}\times\dots\times G_{0}\) for~\(n\ge 1\),
with structure maps
\begin{align}
  \partial_{0}\gg &= [g_{n-2}|\dots|g_{0}], \\
  \partial_{i}\gg &= [\partial_{i-1}g_{n-1}|\dots|\partial_{1}g_{n-i+1}|
  \partial_{0}g_{n-i}\,g_{n-i-1}|\mathrlap{g_{n-i-2}|\dots|g_{0}]} \\
  \notag & & (1\le i\le n-1), \\
  \partial_{n}\gg &= [\partial_{n-1}g_{n-1}|\dots|\partial_{1}g_{1}], \\
  s_{i}\gg &= [s_{i-1}g_{n-1}|\dots|s_{0}g_{n-i}|1_{n-i}|g_{n-i-1}|\dots|g_{0}]
  & (0\le i\le n).
\end{align}
for~\(\gg=[g_{n-1}|\dots|g_{0}]\in(BG)_{n}\).
The canonical twisting function~\(\tau_{G}\colon(BG)_{>0}\to G\) is given by\footnote{%
  The definition in~\cite[p.~88]{May:1968} is incorrect. It only works for right actions on fibres
  whereas the definitions of twisted Cartesian products and twisting functions are based on left actions.}
\begin{equation}
  % \tau_{G}([g_{n-1}|\dots|g_{0}])= g_{n-1}.
  \tau_{G}(\gg)= g_{n-1}^{-1}.
\end{equation}
% for~\([g_{n-1}|\dots|g_{0}]\in(BG)_{n}\).

\subsection{Loop groups}

Let \(X\) be a space that is reduced (meaning that \(X_{0}\) is a singleton).
The \newterm{(Kan) loop group}~\(\Omega X\) of~\(X\) is defined as follows:
\((\Omega X)_{n}\) is the free group generated by~\(X_{n+1}\)
modulo the relations~\(s_{0}x=1_{n}\) for all~\(x\in X_{n}\).
We write \([x]\in(\Omega X)_{n}\) for the equivalence class of~\(x\in X_{n+1}\).
The face and degeneracy maps are determined by
\begin{align}
  \partial_{0}[x] &= [\partial_{0}x]^{-1}[\partial_{1}x], \\
  \partial_{i}[x] &= [\partial_{i+1}x] & (1\le i\le n)\mathrlap{,} \\
  s_{i}[x] &= [s_{i+1}x] & (0\le i\le n)
\end{align}
for~\([x]\in(\Omega X)_{n}\).

By construction, there is a canonical twisting cochain
\begin{equation}
  \tau^{X}\colon X_{>0}\to\Omega X,
  \quad
  x\mapsto[x].
\end{equation}

The classifying space of any group~\(G\) is reduced.
We write \(\bb{g_{n}|\dots|g_{0}}\) for the equivalence class
of~\([g_{n}|\dots|g_{0}]\in(BG)_{n+1}\) in~\((\Omega BG)_{n}\). % (Note that \(BG\) is reduced.)
There is a canonical morphism of groups
\begin{equation}
  \label{eq:OmegaBG-G}
  \Omega BG\to G,
  \quad
  \bb{g_{n}|\dots|g_{0}}\mapsto g_{n}^{-1}
\end{equation}
(which is a homotopy equivalence, compare~\cite[Cor.~27.4]{May:1968}).

\section{Proof of \Cref{thm:main-intro}}
\label{sec:general}

\def\sec#1{\sigma(#1)}
\let\inv\overline

In this section it will be convenient to we write \(\inv{g}\)
instead of~\(g^{-1}\) for~\(g\in G_{n}\).

The projection~\(\pi\colon G\to L\) is a principal \(K\)-bundle.
Let \(\sigma\colon L\to G\) be a pseudo-cross section to~\(\pi\).
We can assume that for any~\(n\ge0\) we have
\begin{equation}
  \label{eq:sigma-1}
  \sigma(1_{n}) = 1_{n}.
\end{equation}
(If \(\rho\colon L\to G\) is any pseudo-cross section,
then we set \(\sigma(l)=\rho(l)\inv{\rho(1_{n})}\) for~\(l\in L_{n}\).)

\begin{proposition}
  \label{thm:action-OmegaBL-BK}
  The loop group~\(\Omega BL\) acts on~\(BK\) via
  \begin{multline*}
    \ll\cdot\kk =
    \Bigl[\, % \partial_{0}\sec{1_{n}}
      \inv{\partial_{0}\sec{l_{n}}}\,k_{n-1}\,\sec{\partial_{0}l_{n}l_{n-1}}\,\inv{\sec{l_{n-1}}}\Bigm| \\
    \partial_{0}\sec{l_{n-1}}\,\inv{\partial_{0}\sec{\partial_{0}l_{n}l_{n-1}}}\,k_{n-2}\,\sec{\partial_{0}^{2}l_{n}\partial_{0}l_{n-1}l_{n-2}}\,\inv{\sec{\partial_{0}l_{n-1}l_{n-2}}}\Bigm|
    \ldots\Bigm| \\
    \partial_{0}\sec{\partial_{0}^{n-2}l_{n-1}\cdots l_{1}}\,\inv{\partial_{0}\sec{\partial_{0}^{n-1}l_{n}\cdots l_{1}}}\,k_{0}\,\sec{\partial_{0}^{n}l_{n}\cdots l_{0}}\,\inv{\sec{\partial_{0}^{n-1}l_{n-1}\cdots l_{0}}}
    \,\Bigr]
  \end{multline*}
  for~\(\ll=\bb{l_{n}|\dots|l_{0}}\in(\Omega BL)_{n}\)
  and~\(\kk=[k_{n-1}|\dots|k_{0}]\in(BK)_{n}\).
\end{proposition}

The arguments to~\(\sigma\) are products of simplices to which
the face operator~\(\partial_{0}\) has been applied repeatedly a decreasing number of times.
For example, \(\partial_{0}^{n-1}l_{n}\cdots l_{1}\) stands
for~\(\partial_{0}^{n-1}l_{n}\,\partial_{0}^{n-2}l_{n-1}\,\cdots\partial_{0}^{2}l_{3}\,\partial_{0}l_{2}\,l_{1}\).
Also note that \(\ll\) acts trivially if~\(l_{n}=1\),
which is necessary as \(\ll=1_{n}\in(\Omega BL)_{n}\) in this case.

\begin{proof}
  This is a direct computation.
\end{proof}

We define the map
\begin{multline}
  \label{eq:def-Psi}
  \Psi\colon BG \to BK \timestauBL BL, \qquad
  \bigl[\,(k_{n-1},l_{n-1})\bigm|\ldots\bigm|(k_{0},l_{0})\,\bigr] \mapsto \\
  \Bigl(
  \bigl[\,k_{n-1}\bigm|
  \partial_{0}\sec{l_{n-1}}\,k_{n-2}\,\sec{l_{n-2}}\,\inv{\sec{\partial_{0}l_{n-1}l_{n-2}}}\bigm| \\
  \partial_{0}\sec{\partial_{0}l_{n-1}l_{n-2}}\,k_{n-3}\,\sec{l_{n-3}}\,\inv{\sec{\partial_{0}^{2}l_{n-1}\partial_{0}l_{n-2}l_{n-3}}}\bigm|
  \ldots\bigm| \\
  \partial_{0}\sec{\partial_{0}^{n-2}l_{n-1}\cdots l_{1}}\,k_{0}\,\sec{l_{0}}\,\inv{\sec{\partial_{0}^{n-1}l_{n-1}\cdots l_{0}}}\,\bigr],
  \bigl[l_{n-1}\bigm|\ldots\bigm|l_{0}\bigr]
  \Bigr).
\end{multline}

\begin{theorem}
  The map~\(\Psi\)
  is an isomorphism of simplicial sets.
  In particular, \(BG\) is a fibre bundle with fibre~\(BK\), base \(BL\)
  and structure group~\(\Omega BL\).
\end{theorem}

\begin{proof}
  It is helpful to use the bijection
  \begin{equation}
    \alpha\colon G\to K\times L,
    \quad
    g \mapsto \bigl(g\,\inv{\sigma\pi(g)}, \pi(g)\bigr).
  \end{equation}
  It commutes with the face maps~\(\partial_{i}\) for~\(i\ge1\) as well as with all
  degeneracy maps.
  Let \(\alpha\in G_{n}\) and~\(g'\in G_{n-1}\).
  If \(\alpha(g)=(k,l)\) and \(\alpha(g')=(k',l')\), then
  \begin{equation}
    \label{eq:alpha-ggprime}
    \alpha(\partial_{0}g\,g') = \bigl(k\,\partial_{0}\sec{l}\,k'\,\sec{l'}\,\inv{\sec{\partial_{0}l l'}} ,l l'\bigr).
  \end{equation}

  The map~\(\alpha\) induces a bijection~\(BG\to BK\times BL\). By what we have just said,
  the transferred structure maps on~\(BK\times BL\) take on the form
  \begin{align}
    \partial_{0}(\kk,\ll) &= (\partial_{0}\kk,\partial_{0}\ll), \\
    \partial_{i}(\kk,\ll) &=
    \Bigl(
    \bigl[
    \partial_{i-1}k_{n-1}\bigm|\ldots\bigm|\partial_{1}k_{n-i+1}\bigm| \\
    \notag & \qquad\quad
    \partial_{0}k_{n-i}\,\partial_{0}\sec{l_{n-i}}\,k_{n-i-1}\,\sec{l_{n-i-1}}\,\mathrlap{\inv{\sec{\partial_{0}l_{n-i}l_{n-i-1}}}\bigm|} \\
    \notag & \qquad\quad
    k_{n-i-2}\bigm|\ldots\bigm|k_{0}
    \bigr]
    ,
    \partial_{i}\ll
    \Bigr)
    & (1\le i\le n-1)\mathrlap{,} \\
    \partial_{n}(\kk,\ll) &= (\partial_{n}\kk,\partial_{n}\ll), \\
    s_{i}(\kk,\ll) &= (s_{i}\kk,s_{i}\ll)& (0\le i\le n)
  \end{align}
  for~\(\kk=[k_{n-1}|\dots|k_{0}]\in BK_{n}\)
  and \(\ll=[l_{n-1}|\dots|l_{0}]\in BL_{n}\).

  It is again a direct computation to verify
  that the map from~\(BK\times BL\) (with the above structure maps)
  to~\(BK \timestauBL BL\) is simplicial.
  The compatibility with degeneracy maps uses \eqref{eq:sigma-1}.  
\end{proof}

\section{Semidirect products}
\label{sec:semidirect}

In this section we assume 
that there is a multiplicative section~\(\sigma\colon L\to G\) to~\(\pi\).
Then \(G\) is isomorphic to the semidirect product~\(K\ltimes L\)
via the map~\((k,l)\mapsto k\,\sigma(l)\).
In this case the formulas given in the preceding section simplify considerably.

We write the conjugation action of~\(L\) on~\(K\)
as \(l*k=\sigma(l)\,k\,\sigma(l)^{-1}\) for~\(k\in K_{n}\) and~\(l\in L_{n}\),
and \(l l'*k\) to mean~\((l l')*k\).
The product in~\(K\ltimes L\) is given by
\begin{equation}
  (k,l)\cdot(k',l') = \bigl(k\,(l*k'),l l'\bigr)
\end{equation}
for~\(k\),~\(k'\in K_{n}\), and \(l\),~\(l'\in L_{n}\).

\begin{corollary}
  \label{thm:action-BK}
  The group~\(L\) acts on~\(BK\) via
  \begin{equation*}
    l \cdot \kk = % \bigl[k_{n-1}\bigm|\cdots\bigm|k_{0}\bigr] =
    \bigl[\,\partial_{0}l*k_{n-1}\bigm|\partial_{0}^{2}l*k_{n-2}\bigm|\ldots\bigm|\partial_{0}^{n}l*k_{0}\,\bigr].
  \end{equation*}
  for~\(l\in L_{n}\) and \(\kk = [k_{n-1}|\dots|k_{0}]\in BK_{n}\).
\end{corollary}

\begin{proof}
  The action of~\(\Omega BL\) on~\(BK\) from \Cref{thm:action-OmegaBL-BK}
  factors through the \(L\) via the canonical morphism~\(\Omega BL\to L\)
  given by~\eqref{eq:OmegaBG-G}.
\end{proof}

Moreover, the map~\(\Psi\) defined in~\eqref{eq:def-Psi} simplifies to
\begin{multline}
  \label{eq:phi-semidirectcase}
  \Phi\colon BG \to BK \timestauL BL, \qquad
  \bigl[\,(k_{n-1},l_{n-1})\bigm|\ldots\bigm|(k_{0},l_{0})\,\bigr] \mapsto \\
  \Bigl(
  \bigl[k_{n-1}\bigm|\partial_{0}l_{n-1}*k_{n-2}\bigm|
  \partial_{0}^{2}l_{n-1}\partial_{0}l_{n-2}*k_{n-3}\bigm|\ldots\bigm|
  \partial_{0}^{n-1}l_{n-1}\cdots\partial_{0}l_{1}*k_{0}\bigr], \\
  \bigl[l_{n-1}\bigm|\ldots\bigm|l_{0}\bigr]
  \Bigr),
\end{multline}
and we have the following:

\begin{corollary}
  \label{thm:main-semidirect}
  The map~\(\Phi\)
  is an isomorphism of simplicial sets.
  In particular, \(BG\) is a fibre bundle with fibre~\(BK\), base \(BL\)
  and structure group~\(L\).
\end{corollary}


\begin{thebibliography}{9}

\bibitem{tomDieck:2008}
T.~tom Dieck,
\newblock \textit{Algebraic topology},
\newblock corrected 2nd printing,
\newblock European Math.\ Soc., Zürich 2010;
\newblock \doi{10.4171/048}

\bibitem{May:1968}
J.\,P.~May,
\newblock \textit{Simplicial objects in algebraic topology},
\newblock Chicago Univ.\ Press, Chicago 1992

\end{thebibliography}
\end{document}